\documentclass[aps,pre,amsmath,amssymb,superscriptaddress,showpacs,twocolumn]{revtex4-1}
\usepackage{graphicx,url}

\begin{document}

\title{Zero forcing number of graphs with a power law degree distribution}

\author{Alexei Vazquez}
\email{avazque1@protonmail.com}
\affiliation{Cancer Research UK Beatson Institute, Glasgow, United Kingdom}
\affiliation{Institute for Cancer Sciences, University of Glasgow, Glasgow, United Kingdom}

\date{\today}

\begin{abstract}
The zero forcing number is the minimum number of black vertices that can turn a white graph black following a single neighbour colour forcing rule. The zero forcing number  provides topological information about linear algebra on graphs, with applications to the controllability of quantum dynamical systems. Here, I investigate the zero forcing number of undirected graphs with a power law degree distribution $p_k\sim k^{-\gamma}$ by means of numerical simulations. For graphs generated by the preferential attachment model, with a diameter scaling logarithmically with the graph size, the zero forcing number approaches the graph size when $\gamma\rightarrow2$. In contrast, for graphs generated by the deactivation model,  with a diameter scaling linearly with the graph size, the zero forcing number is smaller than the graph size independently of $\gamma$.  Therefore the scaling of the graph diameter with the graph size is another factor determining the controllability of dynamical systems. These results have implications for the controllability of quantum dynamics on energy landscapes, often characterized by a complex network of couplings between energy basins. 
\end{abstract}

\maketitle

\section{Introduction}

A zero forcing set of an undirected graph, denoted by $S$, is defined as a subset of vertices that can change a graph coloured white to black according to the following dynamical rules \cite{aim08}. Start with a set of vertices coloured black and all other vertices white. At each step, all black vertices with exactly one white neighbour colour that neighbour black. If the process ends colouring the whole graph black then the starting set of black vertices is a zero forcing set. The zero forcing number of a graph $G$, denoted by $Z(G)$, is defined as the size of the zero forcing set with minimum size.

The zero forcing number has been studied extensively in the context of linear algebra to calculate bounds on the minimum rank of a graph \cite{aim08,hogben10,trefois15}. The minimum rank of a graph $G$, denoted by $mr(G)$, is defined as the  minimum rank of all real matrices with non-zero elements on the graph edges \cite{aim08}. It has been mathematically proven that \cite{hogben10}.
\begin{equation}
mr(G)\geq N - Z 
\label{bound}
\end{equation}
where the equality is warranty when $G$ is an undirected or directed tree. Equation (\ref{bound}) holds true for graphs without loops (simple), with loops, undirected and directed, after small changes in the forcing rule \cite{hogben10}.

More recently, the zero forcing set has found applications in the controllability of quantum dynamical systems \cite{burgarth13}. The linear control theory of quantum systems is based on time-dependent Hamiltonians of the form
\begin{equation}
H(t)=H_I+\sum_k u_k(t)C^{(k)}
\label{Ht}
\end{equation}
operating on a $n$-dimensional Hilbert space ${\cal H}_n$. Here $H_I$ and $C^{(k)}$ are $n\times n$ Hermitian matrices encoding the interactions between internal states and the action of external perturbations, respectively. The Hamiltonian is fully controllable if for every unitary operation $U\in \mathrm{SU}(n)$ there is a control profile $\mathbf{u}(t)$ such that $U={\cal T}\exp[-i\int_0^t H(\tau)d\tau]$, where ${\cal T}$ is the time ordering operator \cite{ramakrishna95}. A sufficient condition for the controllability of the quantum system (\ref{Ht}) is that $H_I$, $C^{(k)}$ and their nested commutators generate the Lie algebra of $\mathrm{SU}(n)$ \cite{ramakrishna95,dalessandro07}. Burgarth {\em et al} \cite{burgarth13} have identified a class of fully controllable quantum systems where the $C^{(k)}$ are associated with the dominating set. The quantum system (\ref{Ht}) has an associated backbone graph $G(H_I)$, where vertices represent the degrees of freedom and edges represent the couplings between them, the non-zero off-diagonal elements of $H_I$. If $S=\{v_1,\ldots,v_m\}$ is a dominating set of $G(H_I)$ then the quantum system (\ref{Ht}) with
\begin{equation}
C^{(k)}_{ij}=\delta_{iv_k}\delta_{jv_k}
\label{CZ}
\end{equation}
is fully controllable \cite{burgarth13}. Obviously, within this class of quantum systems, the minimum number of external inputs needed to fully control the system is the zero forcing number $Z(G)$.

Finding the zero forcing set of undirected graphs is a NP-hard problem \cite{trefois15}. Yet, there can be classes of graphs where heuristic algorithms can find a nearly optimal solution in polynomial time. This is the case for the minimum vertex cover on graphs with a power law degree distribution $p_k\sim k^{-\gamma}$. A vertex cover is a set of coloured vertices in a graph such that every edge is incident to a coloured vertex. The minimum vertex cover is the vertex cover (or covers) of minimum size. Finding the minimum vertex cover is NP-hard as well. Yet, we can use a leaf-removal algorithm to find the minimum vertex cover of graphs with a power law degree distribution \cite{vazquez_weigt03}.

Here I investigate the zero forcing number of graphs with a power law degree distribution $p_k\sim k^{-\gamma}$ with $2<\gamma<\infty$. I introduce a leaf-removal algorithm to identify sub-optimal zero forcing sets. Using graphs generating by preferential attachment and the leaf-removal algorithm, I demonstrate numerically that $Z(G)\rightarrow N$ when $\gamma\rightarrow2$. In contrast, the size of the minimum vertex cover is of order $N(\gamma-2)$ when $\gamma\rightarrow2$.

\begin{figure}[t]
\includegraphics[width=3.3in]{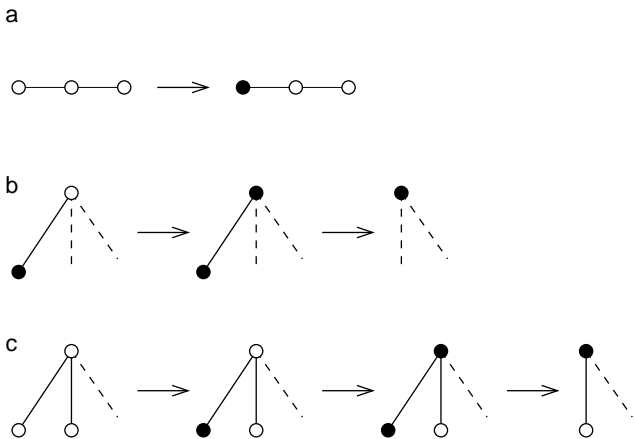}
\caption{Illustration of the leaf removal rules to construct the zero forcing set.}
\label{fig1}
\end{figure}

\section{Leaf removal algorithm}

Leaf-removal algorithms can tackle NP-hard problems on graphs with a power law degree distribution \cite{vazquez_weigt03,zhao15}. A leaf is defined as a vertex with degree 1. The leaf-removal takes advantage of scenarios where the optimal choice is evident. Optimal choice will depend on the specific optimization problem, as shown in Fig. \ref{fig1} for zero forcing. If a leaf is at the end of a chain and the leaf at the other end is white (Fig. \ref{fig1}a), one optimal choice is to colour the leaf black and add the leaf to the forcing set. If a leaf is black and the neighbour is white (Fig. \ref{fig1}b), then the neighbour is forced to black and the leaf is removed. If the neighbour becomes a leaf then the procedure can be continue recursively. Yet another scenario is a vertex with $L$ leafs (Fig. \ref{fig1}c). Here we have no other choice but colouring black $L-1$ leafs and adding them to the zero forcing set. In some instances the leaf removal rules do not create new leafs and the algorithm cannot proceed further. To re-start the process I will remove the vertex with the largest degree and add it to the zero forcing set. Putting all together the algorithm proceed as follows. 
\begin{enumerate}

\item Start with all vertices coloured white and an empty zero forcing set.

\item Create a list of all leafs in the current graph.

\item For each leaf in the list carry on the following rules sequentially.

\begin{itemize}
\item If the leaf is white and it is at the end of an isolated chain with a white leaf at the other end, then colour the leaf black and add it to the zero forcing set.
\item If the leaf is black, remove the leaf and colour the neighbour black. In the event that the neighbour becomes itself a leaf proceed recursively until no new leaf is created (Fig. \ref{fig1}a).
\item Otherwise, find all leafs at distance 2 from the current leaf, colour them black, remove them and add them to the zero forcing set. If no leaf is found at distance 2 no action is taken.
\end{itemize}

\item If the graph is not empty, find a vertex with the current largest degree, colour the vertex black, remove the vertex, add the vertex to the zero forcing set, and move all neighbour vertices that become a leaf to the leaf list.

\item Remove all isolated vertices, those that are white colour them black and add them to the zero forcing set.

\item If the graph is empty stop, otherwise go back to step 2.

\end{enumerate}
The size of the resulting zero forcing set will be denoted by $Z_{LM}(G)$, where $LM$ stands for leaf and maximum degree removal. Since the maximum degree removal, step 4, is not necessarily optimal, this algorithm overestimates the zero forcing number,
\begin{equation}
Z_{LM} \geq Z
\label{ZLMbound}
\end{equation}

As a comparison, I will also calculate the minimum vertex covering using an adaptation of the vertex covering leaf removal algorithm \cite{vazquez_weigt03}. The algorithm proceed as follows. Start with all vertices uncovered and a list of leafs in the graph. If the leaf list is non-empty, extract a leaf, remove the leaf, add the leaf neighbour to the vertex covering set and remove the leaf neighbour. Otherwise, find a vertex with the current largest degree, add it to the vertex covering set and remove it. Continue until the graph is empty. The minimum vertex covering will be denoted by $V(G)$. The vertex cover size estimated by the leaf-maximum degree removal will be denoted by $V_{LM}(G)$.  Once again, because of the maximum degree removal rule, this algorithm overestimates the minimum vertex covering,
\begin{equation}
V_{LM} \geq V
\label{ZLMbound}
\end{equation}

\begin{figure}[t!]
\includegraphics[width=3.3in]{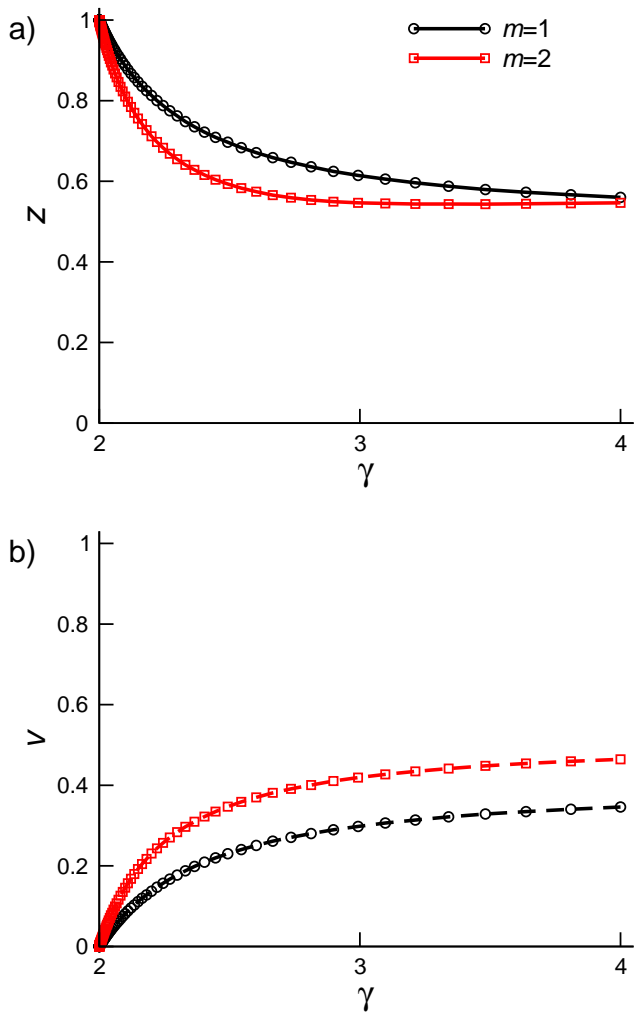}
\caption{Zero forcing fraction $z=Z_{LM}/N$ and minimum vertex cover fraction $v=V_LM/N$ of preferential attachment graphs, with $N=10,000$ and averaged over 100 graphs.}
\label{fig2}
\end{figure}

\section{Preferential-attachment graphs}

First, I will consider graphs generated by the preferential attachment model with initial attractiveness \cite{dorogovtsev00} ($G_{PA}$). The graph is started with a fully connected graph of $m+1$ vertices. Then add new vertices one by one until the targeted graph size is reached. Each time a new vertex is added it is connected to $m$ existing and non-overlapping vertices, each selected with probability
\begin{equation}
\pi_i = \frac{ (a-1) m + k_i }{ \sum_j \left[ (a-1) m + k_j \right] }
\label{preferential_attachment}
\end{equation}
where the indexes run over vertices in the current graph and $k_i$ denotes the degree of vertex $i$. The parameter $a>0$ represents a vertex independent attractiveness named initial attractiveness \cite{dorogovtsev00}. This model generates graphs with a power law degree distribution $p_k \sim k^{-\gamma}$ with exponent \cite{dorogovtsev00}
\begin{equation}
\gamma = 2 + a
\label{gammaPA}
\end{equation}
By tuning the initial attractiveness $a$ we can obtain exponents in the range $2< \gamma < \infty$. The case $a=1$ ($\gamma=3$) corresponds to the original Barab\'asi-Albert model \cite{barabasi99}.

\begin{figure}[t!]
\includegraphics[width=3.3in]{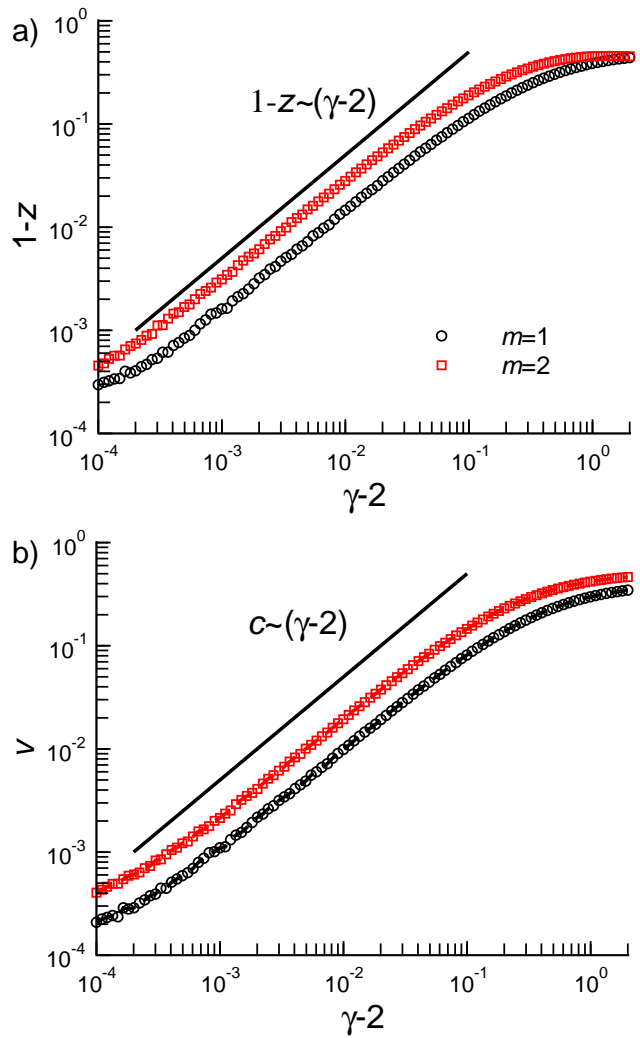}
\caption{Scaling of the zero forcing fraction $z=Z_{LM}/N$) and minimum vertex cover fraction $v=V_LM/N$ of preferential attachment graphs, with $N=10,000$ and averaged over 100 graphs.}
\label{fig3}
\end{figure}

Using the LM algorithms, I have estimated the zero forcing number and minimal vertex covering of graphs generated by the preferential attachment model (Fig. \ref{fig2}). The zero forcing number and the minimum vertex covering behave quite differently as $\gamma\rightarrow2$.  The zero forcing number approaches the graph size while the minimum vertex covering approaches zero. Around $\gamma=2$ we observe the scalings (Fig. \ref{fig3})
\begin{equation}
N-Z_{LM}(G_{PA}) \approx c_z N (\gamma-2)
\label{Zscaling}
\end{equation}
\begin{equation}
V_{LM}(G_{PA}) \approx c_v N (\gamma-2)
\label{Vscaling}
\end{equation}
Substituting the scaling (\ref{Zscaling}) into the minimum rank lower bound (\ref{bound}) we obtain
\begin{equation}
mr(G_{PA}) \geq c_z N (\gamma-2)
\label{mrscaling}
\end{equation}
When $\gamma\approx2$ we cannot exclude that the minimum rank is 0.

\begin{figure}[t]
\includegraphics[width=3.3in]{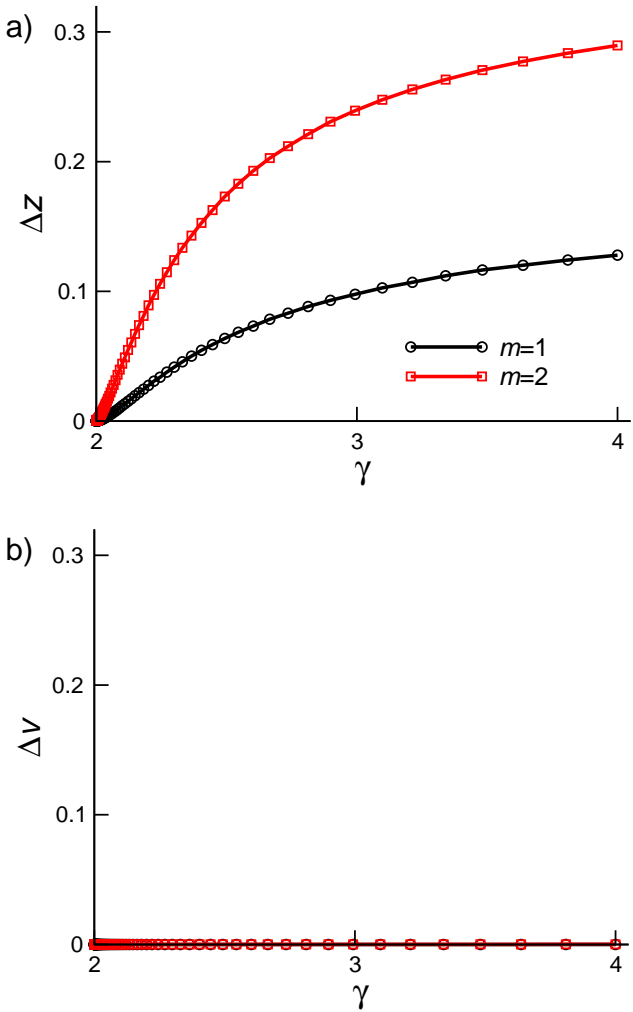}
\caption{Contribution of the maximum degree step to the size of the zero forcing and vertex covering sets estimated by the LM algorithm, for preferential attachment graphs with $N=10,000$ and averaged over 100 graphs.}
\label{fig4}
\end{figure}

Since the LM algorithm yields upper bounds it is worth asking how tight are those bounds, specially for $\gamma\approx2$. To address this question I have calculated the number of vertices that were forced using the maximum degree step, denoted by $\Delta Z(G)$ and $\Delta V(G)$ for the zero forcing and vertex covering algorithms, respectively. We observe that $\Delta V(G_{PA})/N\approx0$ for all values of $\gamma$ explored (Fig. \ref{fig4}b). Therefore
\begin{equation}
V_{LM}(G_{PA}) \approx V(G_{PA})
\label{VV_PA}
\end{equation}
In contrast,  $\Delta Z_{LM}(G_{PA})>0$ for most values of $\gamma$ and we cannot tell how good are the zero forcing number estimates. Nevertheless, for $\gamma\rightarrow2$ we have  $\Delta Z_{LM}(G_{PA})/N\rightarrow0$. In the vicinity of $\gamma=2$ the LM algorithm provides good estimates of the zero forcing number, {\em i.e.} 
\begin{equation}
Z_{LM}(G_{PA}) \approx Z(G_{PA})
\label{VV_PA}
\end{equation}
for $\gamma\approx2$. This allows us to conclude that the zero forcing number of preferential-attachment graphs approaches the graph size when $\gamma\rightarrow2$.

\begin{figure}[t]
\includegraphics[width=3.3in]{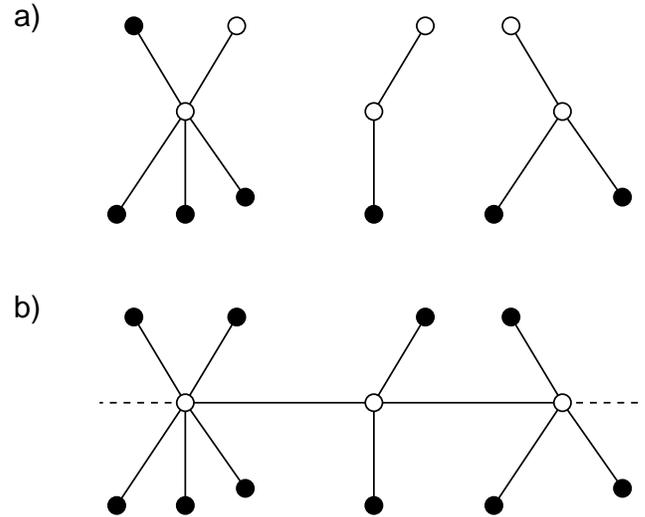}
\caption{Star graphs. a) Isolated stars. b) String of stars. Solid circles are vertices in the zero forcing set.}
\label{fig_star}
\end{figure}

\section{Star and deactivation graphs}

It is straightforward to show that $\gamma\rightarrow2$ is not a sufficient condition to obtain $Z(G)\approx N(G)$. To this end, consider graphs made by isolated starts with a degree distribution $p_k$, $k=1,\ldots$ ($G_{0\star}$, Fig. \ref{fig_star}a). In this case the degree distribution satisfy the constraint
\begin{equation}
\sum_{k>1} p_k k = p_1
\label{stars_dis}
\end{equation}
For each isolated star all but 1 leaf need to be included in the zero forcing set, resulting in
\begin{equation}
\frac{Z(G_{0\star})}{N} = \sum_{k>1} p_k (k-1) =  2p_1-1 < 1
\label{Zstar}
\end{equation}

A similar result is obtained for graphs made by a string of stars ($G{-\star}$, Fig. \ref{fig_star}b). In this case the degree distribution satisfy the constraint
\begin{equation}
\sum_{k>1} p_k (k-2) = p_1
\label{stars_dis_1}
\end{equation}
In this case all leafs need to be forced resulting in  
\begin{equation}
\frac{Z(G_{-\star})}{N} = \sum_{k>1} p_k (k-2) =  p_1 < 1
\label{Zstar}
\end{equation}
$Z(G_{0\star})<N$ and $Z(G_{-\star})<N$ regardless of the shape of the degree distribution. A power low degree distribution $p_k\sim k^{-\gamma}$ with $\gamma\rightarrow2$ is not a sufficient condition for $Z(G)=N(G)$.

\begin{figure}[t]
\includegraphics[width=3.3in]{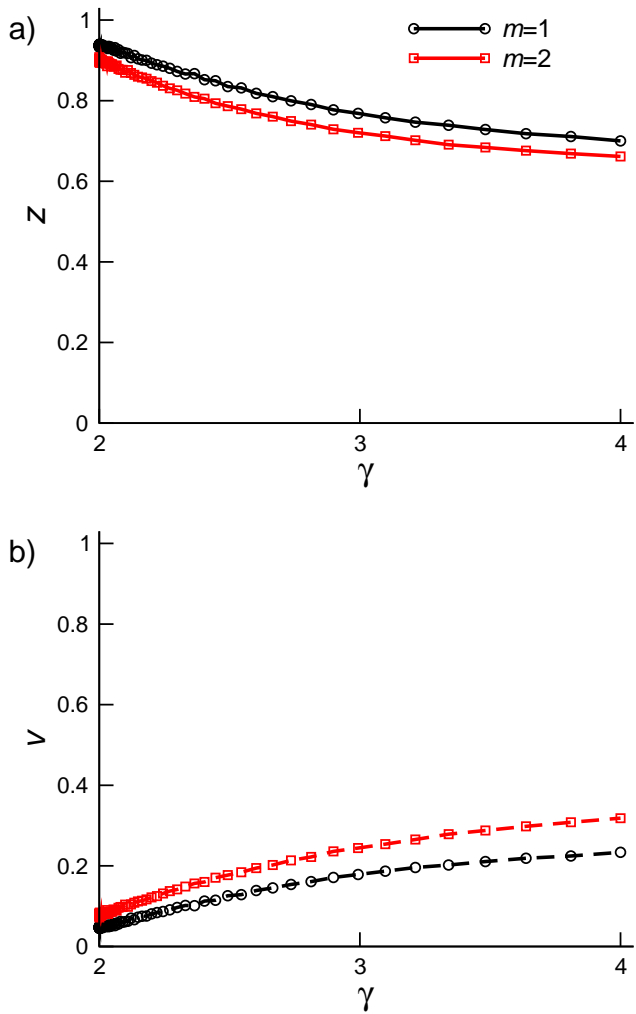}
\caption{Zero forcing fraction $z=Z_{LM}/N$ and minimum vertex cover fraction $v=V_LM/N$ of deactivation graphs, with $N=10,000$ and averaged over 100 graphs.}
\label{fig5}
\end{figure}

The deactivation graphs of Klemm and Egu\'{\i}luz \cite{klemm02} ($G_D$) provide another counter example inspired on natural rules of network evolution. The basic idea of the deativation graph model is that, with time, some vertices will no longer participate in the network evolution, becoming inactive or deactivated. That rule models the retirement of a scientist in the context of co-authorship networks for example. The deactivation graphs studied here are generated as follows. Start with a fully connected graph of $m+1$ active vertices. At each graph evolution step, add a new active vertex, connect the new vertex to the pre-existing $m$ active vertices and set one of the active vertices ($i\in{\cal A}$) inactive with a probability
\begin{equation}
\pi_i = \frac{ \sum_{s\in{\cal A}} ((a-1)m+k_s)^{-1} }{ (a-1)m + k_i }
\label{d}
\end{equation}
The deactivation model generate graphs with a power low degree distribution $p_k\sim k^{-\gamma}$ with exponent \cite{klemm02}
\begin{equation}
\gamma = 2 + a
\label{gamma_d}
\end{equation}
By tuning $a$ we can thus generate power law exponents in the range $2\leq\gamma<\infty$.

Using the leaf-maximum degree removal algorithms, I have estimated the zero forcing number and minimal vertex covering of deactivation graphs (Fig. \ref{fig5}). For the deactivation graphs the zero forcing number does not approach the graph size when $\gamma\rightarrow2$ (Fig. \ref{fig5}a). In fact, $Z_{LM}(G_D)/N<1$ and $V_{LM}(G_D)/N>0$ for all values of $\gamma$ (Fig. \ref{fig5}a,b).

As shown before, there is a key difference between the preferential-attachment and deactivation graphs regrading the graph diameter, denoted by $d$. The preferential-attachment generate small-world graphs for $\gamma>3$ ($d\sim\ln N$) and ultra-small graphs ($d\sim\ln\ln N$) for $2<\gamma<3$ \cite{cohen03}. In contrast, deactivation graphs are effectively one-dimensional ($d\sim N$) \cite{vazquez_deactivation03}. This, together with the analysis of star graphs, indicates that the small-world property is a requirement to obtain $Z(G)\approx N$ when $\gamma\rightarrow2$.

\section{Conclusions}

Using numerical simulations of preferential-attachment graphs, I have demonstrated that the zero forcing number approaches the graph size ($Z\rightarrow N$) when the exponent of the power law degree distribution approaches 2 ($\gamma\rightarrow2)$. This extends a similar result by Liu-Slotine-Barab\'asi for directed graphs \cite{liu11} to undirected graphs. Through the analysis of some counterexamples, I have shown that the small-world property of preferential-attachment graphs is a necessary requirement for this result.

These observations are relevant for the controllability of quantum systems. If the degrees of freedom are coupled by a small world network with a power law degree distribution, then we would need to manipulate almost all degree of freedoms to achieve full controllability. A putative scenario is the configuration space of energy basins associated with protein folding. The network of energy basins of a Lennard-Jones gas and a $\beta$-sheet peptide is characterized by a power law degree distribution,  with exponents $\gamma\approx2.8$ \cite{doye02} and $\gamma\approx2$ \cite{rao04}, respectively. For the Lenard-Jones gas it was further shown that the network is small world. Controlling the quantum dynamical evolution of these systems with the zero forcing Hamiltonians (\ref{CZ}) would be unfeasible. It would require the manipulation of almost all energy basins to achieve full control. 

\bibliographystyle{apsrev4-1}


\bibliography{forcing}

\end{document}